\documentclass[12pt,amsart]{amsart}
\usepackage{amsmath,amsxtra,amstext,amssymb,latexsym,amsthm,dsfont,amscd,xypic}
\usepackage{setspace}
\usepackage{cite}
\newtheorem{thm}{Theorem}
\newtheorem{lemma}[thm]{Lemma}
\newtheorem{cor}[thm]{Corollary}
\newtheorem{prop}[thm]{Proposition}
\newtheorem{definition}[thm]{Definition}
\newtheorem{rem}[thm]{Remark}

\newcommand{\Tv}{T_{V_{\bullet}}}
\newcommand{\Fla}{{\rm F\ell}(a,n)}
\newcommand{\Fln}{{\rm F\ell}(n)}
\newcommand{\Fl}{{\rm F\ell}}
\newcommand{\bull}{_{\bullet}}
\newcommand{\SLn}{SL_n}
\newcommand{\tilE}{\tilde{E}_{\bullet}}

\newcommand{\codim}{{\rm codim}}

\newcommand{\Hom}{{\rm Hom}}
\newcommand{\Gr}{{\rm Gr}}

\begin{document}

\title[Horn recursion for $H^*(\Fla)$]{A partial Horn recursion in the cohomology of flag varieties}

\author{Edward Richmond}
\address{Department of Mathematics\\ UNC-Chapel Hill\\ CB \#3250, Phillips Hall \\ Chapel Hill,  NC 27599}
\email{erichmon@email.unc.edu}

\maketitle

\begin{abstract}Horn recursion is a term used to describe when non-vanishing products of Schubert classes in the cohomology of complex flag varieties are characterized by inequalities parameterized by similar non-vanishing products in the cohomology of ``smaller" flag varieties.  We consider the type A partial flag variety and find that its cohomology exhibits a Horn recursion on a certain deformation of the cup product defined by Belkale and Kumar in ~\cite{BK06}.  We also show that if a product of Schubert classes is non-vanishing on this deformation, then the associated structure constant can be written in terms of structure constants coming from induced Grassmannians.\end{abstract}


\section{Introduction}\label{SECTION_INTRO}

\quad One of the primary concerns of Schubert calculus is to determine the product structure of the cohomology ring of flag varieties with respect to its basis of Schubert classes.  Let $P$ be a parabolic subgroup of some semi-simple complex connected algebraic group $G$ and consider the flag variety $G/P$.  Suppose we are given Schubert varieties $X_{w^1},X_{w^2},\ldots, X_{w^s}\subseteq G/P$ such that the product \begin{equation}\label{equation_product}\prod_{k=1}^s [X_{w^k}]=c[pt]\in H^*(G/P)=H^*(G/P,\mathds{Z})\end{equation} for some integer structure coefficient $c\in\mathds{Z}_{\geq 0}$.  The two questions we ask are: under what conditions is $c\neq 0$ and more specifically, can we explicitly compute $c$?  Horn's conjecture~\cite{Ho62} on the Hermitian eigenvalue problem provides an answer to the first question in the case of the usual Grassmannian.  The result is that $c\neq 0$ if and only if the Schubert data $(w_1,\ldots,w_s)$ satisfies a certain list of linear inequalities called Horn's inequalities. Remarkably, Horn's inequalities themselves are indexed by the non-vanishing of similar structure coefficients for smaller Grassmannians (Horn recursion). Horn's conjecture was proved by the work of Klyachko~\cite{Kly98} and the saturation theorem of Knutson and Tao~\cite{KT99}. Belkale~\cite{Be06} later gives a geometric proof of Horn's conjecture set in the context of intersection theory.  Purbhoo and Sottile~\cite{PS06} have shown that any cominuscule flag variety exhibits a Horn recursion as well.  However, Horn recursion has not yet been determined in any non-cominuscule flag varieties.  A general discussion about the connection of Horn's inequalities to other topics can be found in Fulton~\cite{Fu00, Fu00a}.
\smallskip



\quad In this paper we consider the case where $G=\SLn$ and $P$ is any parabolic subgroup.   The homogeneous space $\SLn/P$ corresponds to a partial flag variety $\Fla$ for some sequence of integers $a=\{0<a_1<a_2<\cdots<a_r<n\}$ and the Schubert classes are indexed by the set $$S_n(a):=\{(w(1),w(2),\ldots,w(n))\in S_n \ |\ w(i)<w(i+1) \ \forall \ i\notin \{a_1,a_2,\ldots,a_r\}\}$$  where $S_n$ denotes the permutation group on $[n]:=\{1,2,\ldots,n\}$.  In the case of the Grassmannian, we will use the notation $\Gr(r,n)$ and denote the set $S_n(\{r\})$ by simply $S_n(r)$.  We find that $\Fla$ exhibits a Horn recursion on a certain deformation of the cup product and that the structure coefficients on this deformation are a product of structure coefficients coming from induced Grassmannians.
\smallskip



\subsection{A Horn recursion for partial flag varieties}
For any $s$-tuple $(w^1,w^2,\ldots,w^s)\in S_n(a)^s$ such that $\sum^s_{k=1}\codim\, X_{w^k}=\dim \Fla$, we define the associated structure constant to be the integer $c$ given in equation \eqref{equation_product}.  Assume that the structure constant associated to $(w^1,w^2,\ldots,w^s)$ is nonzero.  For any $i\in [r]$, consider the projection $f_i:\Fla \rightarrow \Gr(a_i,n)$.  The expected dimension of the intersection $\bigcap_{k=1}^s f_i(X_{w^k})$ in $\Gr(a_i,n)$ is nonnegative and hence we have
\begin{equation}\label{numerical_cond_ineq}
\sum_{k=1}^s\sum_{j=1}^{a_i} \Big(n-a_i+j-w^k(j)\Big) \leq a_i(n-a_i)\end{equation}

for all $i\in [r]$.  Note that the $LHS$ of equation \eqref{numerical_cond_ineq} is equal to the sum of the codimensions of the varieties $f_i(X_{w^k})\subseteq\Gr(a_i,n)$.  We find that $\Fla$ exhibits a Horn recursion on the boundary where the inequalities \eqref{numerical_cond_ineq} are equalities. In other words, where the expected dimension of $\bigcap_{k=1}^s f_i(X_{w^k})$ is zero for all $i\in\, [r]$. Equivalently, this numerical condition can be described by the notion of \emph{Levi-movability} defined by Belkale and Kumar in ~\cite{BK06}.  Let $\tilE:=E_{a_1}\subseteq \cdots \subseteq E_{a_r}\in\Fla$ denote the standard partial flag which is identified with the point $eP\in\SLn/P\simeq\Fla$ and let $L$ be the Levi subgroup of $P\subseteq \SLn$. We say $(w^1,w^2,\ldots,w^s)\in S_n(a)^s$ is Levi-movable (or $L$-movable) if for generic $(l_1,l_2,\ldots,l_s)\in L^s$ the intersection $\bigcap_{i=1}^s l_i(w^i)^{-1}X_{w^i}$ is finite and transverse at $\tilE\in\Fla$.  This condition on $(w^1,w^2,\ldots,w^s)$ can be quite restrictive.  For instance, if $\Fla$ is the complete flag variety, then any structure coefficient associated to a $L$-movable $s$-tuple is equal to 1 (this result is stated in Corollary \ref{Cor_Stcnst1} of Theorem \ref{Thm_constants}).   The notion of $L$-movability defines a new product on the cohomology of $\Fla$ which is a deformation of the usual product.  For more details on this new product see ~\cite{BK06}.  Before we state the first main result of this paper, we need the following definition of induced permutations:

\begin{definition}\label{def_induced_Weyl}Let $A=\{\dot{a}_1< \dot{a}_2< \cdots< \dot{a}_d\}\subseteq [n]$ be an ordered subset of cardinality $d<n$ and $w\in S_n$.  Define $w_A\in S_d$ to be the unique permutation with the same descent set as the sequence $\{w(\dot{a}_1), w(\dot{a}_2), \ldots, w(\dot{a}_d)\}$.\end{definition}


For example if $w=(2,5,3,1,4)\in S_5$ and $A=\{1,2,5\}\subseteq [5]$, then $$w=(\, \boxed{2,5,}\, 3,1,\boxed{4}\,)\mapsto w_A=(1,3,2)\in S_3.$$  The map $w\mapsto w_A$ is the same as the flattening function of Billey and Braden in ~\cite{BB03}.  Set $a_0=0$ and $a_{r+1}=n$ and define $b_i:=a_i-a_{i-1}$ and $A_i:=\{a_{i-1}+1<a_{i-1}+2< \cdots<a_i\}\subset [n]$.  For any $w\in S_n(a)$ and $\{i<j\}\subseteq [r+1]$, define $w_{i,j}:=w_{A_i\cup A_j}\in S_{b_i+b_j}$.

\begin{thm}\label{Thm1} Let $(w^1, w^2,\ldots , w^s) \in S_n(a)^s$ such that $$\sum^s_{k=1}\codim\, X_{w^k}=\dim\Fla.$$
Then the following are equivalent:
\begin{enumerate}\item[(i)] $\prod_{k=1}^s [X_{w^k}]$ is a nonzero multiple of a class of the point in $H^*(\Fla)$ and
\begin{equation}\label{numerical_cond_eq}\sum_{k=1}^s\sum_{j=1}^{a_i} \Big(n-a_i+j-w^k(j)\Big) = a_i(n-a_i)\quad \forall\ i\in\, [r].\end{equation}
\item[(ii)] The s-tuple $(w^1, w^2,\ldots , w^s)$ is $L$-movable.
\item[(iii)] $\prod_{k=1}^s [X_{w^k_{i,j}}]$ is a nonzero multiple of a class of the point in $H^*(\Gr(b_i, b_i+b_j))\quad \forall\ \{i<j\}\subseteq [r+1]$.
\item[(iv)] For any $\{i<j\}\subseteq[r+1]$, the following are true:
 \begin{enumerate}\item The sum $\sum_{k=1}^s\codim(X_{w^k_{i,j}})=\dim\Gr(b_i,b_i+b_j)$.\item For any $1 \leq d < b_i$ and any choice $(u^1, u^2,\ldots , u^s) \in S_{d}(b_i)^s$ such that $\prod_{k=1}^s [X_{u^k}]$ is a nonzero multiple of a class of the point in $H^*(\Gr(d, b_i))$, the following inequality is valid:$$ \sum_{k=1}^s \sum_{l=1}^{d} \Big(b_j+u^k(l)-w^k_{i,j}(u^k(l))\Big) \leq db_j.$$\end{enumerate}\end{enumerate}\end{thm}

Note that $(iii)\Leftrightarrow (iv)$ is a direct application of Horn's conjecture ~\cite{Kly98, KT99} applied to the Grassmannians $\Gr(b_i, b_i+b_j)$ and that $(i)\Leftrightarrow (ii)$ is proved in ~\cite[Thm 15]{BK06}.  We also remark that $(ii)\Rightarrow(iv)$ was proved in ~\cite[Thm 32]{BK06} using geometric invariant theory.   In this paper we will prove Theorem \ref{Thm1} by focusing on $(ii)\Leftrightarrow (iii)$, however only $(iii)\Rightarrow(ii)$ is not previously known.
\smallskip


\quad The main object we use to prove Theorem \ref{Thm1} is the tangent space $T_{\tilE}\Fla$. We find that the subtangent spaces of Schubert cells corresponding to a $L$-movable $s$-tuple decompose into a direct sum with respect to a certain decomposition of the tangent space of $\Fla$. Fix a splitting $Q_1\oplus Q_2 \oplus \cdots \oplus Q_{r+1}$ of $\mathds{C}^n$ such that $E_{a_k}=Q_1\oplus Q_2 \oplus \cdots \oplus Q_{k}.$  Then there exists a natural decomposition $$T_{\tilE}\Fla \ \simeq \ \bigoplus_{i<j}^{r+1} \Hom(Q_i, Q_j).$$  In Proposition \ref{Prop_Tan_decomp}, we show that for any $l\in L$ and $w\in S_n(a)$, $$T_{\tilE}(lw^{-1}X_w)\ \simeq \ \bigoplus_{i<j}^{r+1} \big(\Hom(Q_i,Q_j)\cap T_{\tilde{E\bull}}(lw^{-1}X_w)\big).$$  Hence $(w^1,w^2,\ldots,w^s)$ is $L$-movable if and only if $\bigcap_{i=1}^s T_{\tilE}(l_i(w^i)^{-1}X_{w^i})$ is transversal in each summand $\Hom(Q_i,Q_j) \subseteq T_{\tilE}\Fla$ for generic $(l_1,l_2,\ldots,l_s)\in L^s$. We then show that for each summand we can reduce to the case of a Grassmannian where sufficiency is already known to be true.

\subsection{Computing structure coefficients}
Let $c^v_{w,u}$ be the structure coefficients defined by the product $$[X_w]\cdot[X_u]=\sum_{v\in S_n(a)}c^v_{w,u}[X_v].$$   In the case of the Grassmannian, the coefficients $c^v_{w,u}$ are the Littlewood-Richardson numbers which have several nice combinatorial formulas.  However computing these coefficients for $\Fla$ is much more difficult.  We find that if $c^v_{w,u}$ comes from an $L$-movable triple, then $c^v_{w,u}$ is a product of Littlewood-Richardson numbers.  It is well known that the coefficient $c^v_{w,u}$ is the number of points in a generic intersection of the associated Schubert varieties if this intersection is finite. Consider the projection of $f_1:\Fla\twoheadrightarrow\Gr(a_1,n)$. For any Schubert variety $X_w\subseteq\Fla$ we have the induced Schubert varieties $X_{w_1}:=f_1(X_w) \subseteq \Gr(a_1,n)$ and $X_{w_{\gamma}}:=X_w\cap f_1^{-1}(V)\subseteq f_1^{-1}(V)$ in the fiber over a generic point $V\in X_{w_1}$. It is easy to see that $f_1^{-1}(V)\simeq {\rm{F}}\ell(a_\gamma,n-a_1)$ where $a_\gamma=\{0<a_2-a_1<\cdots<a_r-a_1<n-a_1\}$.  We now state the second main result of this paper.

\begin{thm}\label{Thm_constants} Let $(w^1,w^2,\ldots,w^s)\in S_n(a)^s$ be $L$-movable. If $c,c_1$ and $c_\gamma$ are the structure coefficients:
$$\prod_{k=1}^s [X_{w^k}]=c[X_e],\qquad \prod_{k=1}^s [X_{w_1^k}]=c_1[X_e],\qquad \prod_{k=1}^s [X_{w_{\gamma}^k}]=c_\gamma[X_e],$$
in $H^*(\Fla), H^*(\Gr(a_1,n))$ and $H^*(\Fl(a_\gamma,n-a_1))$ respectively, then $c=c_1\cdot c_\gamma.$\end{thm}

In Proposition \ref{L_induced}, we show that the induced $s$-tuple $(w_{\gamma}^1,w_{\gamma}^2,\ldots,w_{\gamma}^s)$ is also Levi-movable.  Hence we can once again use Theorem \ref{Thm_constants} to decompose $c_{\gamma}$. Thus finding the $L$-movable structure coefficients completely reduces to the case of finding these coefficients in the Grassmannians $\Gr(b_i,n-a_{i-1})$ for $i\in[r]$. Note that these are not the same Grassmannians found in Theorem \ref{Thm1}$(iii)$.  In the case of the complete flag variety, we have that $b_i=a_i-a_{i-1}=1$ for all $i\in[n]$. Hence $\Gr(b_i,n-a_{i-1})$ is a projective space where all structure coefficients are equal to 1.  Applying Theorem \ref{Thm_constants} yields the following corollary.

\begin{cor}\label{Cor_Stcnst1} If $\Fla$ is the complete flag variety, then any structure coefficient associated to a $L$-movable $s$-tuple is equal to 1.\end{cor}


\section{Preliminaries on partial flag varieties}\label{FLASL}

\quad Let $\SLn$ be the special linear group acting on the vector space $\mathds{C}^n$.  Let $H \subseteq \SLn$ be the standard maximal torus of diagonal matrices of determinant 1 and let $B$ denote the standard Borel subgroup of upper triangular matrices (of determinant 1) containing $H$.  It is well known that the Weyl group $W:=N(H)/H$ can be identified with $S_n$ the permutation group on $n$ elements. Let $P$ be a parabolic subgroup containing $B$ and let $L=L_P$ be the Levi (maximal reductive) subgroup of $P$.  Let $R$ be the set of roots and let $R^+$ denote the set of positive roots.  Let $\Delta = \{\alpha_1,\alpha_2,\ldots\alpha_{n-1}\}\subset R^+$ be the set of simple roots.

\smallskip

\quad Let $\Delta(P)\subset\Delta$ denote the simple roots associated to $P$ and let $\Delta\backslash\Delta(P)=\{\alpha_{a_1}, \alpha_{a_2}, \ldots, \alpha_{a_r}\}$.  We associate $P$ to the subset $\{a_1, a_2, \ldots, a_r\}\subseteq [n-1]$.  Let $W_P\subseteq W$ be the Weyl group of $P$ generated by the simple transpositions $\{(i, i+1)\ |\ \alpha_i \in \Delta(P)\}$. Define $W^P$ to be the set of minimal length representatives of the coset space $W/W_P$ where the length of $w$ is defined as $\ell(w) := \#\{i<j \ | \ w(i)>w(j)\}.$  We will denote $W^P$ by $S_n(a)$ defined as $$S_n(a):=\{(w(1),w(2),\ldots,w(n))\in S_n \ |\ w(i)<w(i+1) \ \forall \ i\notin \{a_1,a_2,\ldots,a_r\}\}.$$ It is easy to see that $W^P=S_n(a)$ as subsets of $S_n$.  Consider the homogeneous space $\SLn/P$.  This space is $\SLn$-equivariantly isomorphic to the variety of partial flags $$\Fla := \{ V\bull = V_1 \subset V_2 \subset \cdots \subset V_r \subset\mathds{C}^n\ |\ \dim(V_i)=a_i \}$$  The dimension of $\Fla$ is equal to $\sum_{i=1}^r a_i(a_{i+1}-a_i)$ where we set $a_0=0$ and $a_{r+1}=n$. If $a=[n-1]$, we will denote the variety of complete flags by $\Fln$.  In general, if $V$ is a complex vector space let $\Fl(V)$ denote the variety of complete flags on $V$. For any complete flag $F\bull = F_1\subset F_2\subset \ldots \subset F_n\in \Fln$ and $w \in S_n(a)$, we define the Schubert cell by $$X_w^{\circ}(F\bull):= \{ V\bull \in \Fla \ | \dim(V_i \cap F_j) = \# \{t\leq a_i : w(t) \leq j \} \quad \forall i,j \}.$$ Define the Schubert variety to be the closure $X_w(F\bull):=\overline{X_w^{\circ}(F\bull)}\subseteq\Fla$.  The dimension of $X_w(F\bull)$ is equal to $\ell(w)$.  We denote the cohomology class of $X_w(F\bull) \subseteq \Fla$ by $[X_{w}]$. Note that $[X_{w}]$ is independent of the choice of flag $F\bull\in\Fln$ and that the set $\{[X_{w}]\ |\ w\in S_n(a)\}$ forms an additive basis for $H^*(\Fla)$.

\subsection{\textit{Transversal Intersections}}
The condition of a transversal intersection of varieties is one on the tangent spaces. Let $X_1, X_2, \ldots, X_s$ be smooth connected subvarieties of smooth variety $X$.  We say an intersection is transversal at a point $x\in \bigcap_{k=1}^s X_k$, if $$\dim\Big(\bigcap_{k=1}^s T_xX_k\Big)=\dim(T_xX)-\sum_{k=1}^s \codim(T_xX_k).$$ The following proposition is a basic fact on transversality and Schubert varieties.  Let $\mathcal{F}:=(F\bull^1,F\bull^2,\ldots,F\bull^s)$ denote a point in $\Fln^s$.

\begin{prop}\label{Prop Ful_Mac}
Let $(w^1,w^2,\ldots,w^s) \in S_n(a)^s$.  The following are equivalent.
\begin{enumerate}
\item[(i)] $\prod_{k=1}^s [X_{w^k}]\neq0$ in $H^*(\Fla)$.
\item[(ii)] There exists $V\bull\in\Fla$ and $\mathcal{F}\in\Fln^s$ such that $V\bull\in\bigcap_{k=1}^s X_{w^k}^{\circ}(F^k\bull)$ and the Schubert cells $X_{w^k}^{\circ}(F^k\bull)$ meet transversally at $V\bull$.
\item[(iii)] For any $V\bull \in \Fla$, there exists $\mathcal{F}\in \Fln^s$ such that $V\bull\in \bigcap_{k=1}^s X_{w^k}^{\circ}(F^k\bull)$ and the Schubert cells $X_{w^k}^{\circ}(F^k\bull)$ meet transversally at $V\bull$.\end{enumerate}\end{prop}

Note that $\SLn$ acts transitively on $\Fla$ and hence part $(ii)$ is trivially equivalent to part $(iii)$ in Proposition \ref{Prop Ful_Mac}.

\subsection{\textit{Levi-movability}}

Let $\{e_1,e_2,\ldots,e_n\}$ be the standard basis for $\mathds{C}^n$.  For any $i\in [r+1]$, the action of $L$ fixes the subspace $Q_i$ where $Q_i:=\langle e_{a_{i-1}+1}, e_{a_{i-1}+2}, \ldots, e_{a_i}\rangle$. This defines a splitting of

\begin{equation}\label{Q_splitting}\mathds{C}^n=Q_1\oplus Q_2\oplus \cdots \oplus Q_{r+1}.\end{equation}

It is easy to see that the Levi subgroup decomposes into the product
\begin{equation}\label{Levi_split} L\simeq\{(g_1,\ldots,g_{r+1})\in GL(Q_1)\times\cdots\times GL(Q_{r+1}) \ |\  \prod_{i=1}^{r+1}\det(g_i)=1 \}.\end{equation}

Consider the standard complete flag $E\bull \in\Fln$ defined by $E_i:=\langle e_1, e_2, \ldots e_i\rangle$ and partial sub-flag $$\tilE:=E_{a_1}\subseteq E_{a_2}\subseteq\cdots\subseteq E_{a_r}\in\Fla.$$  For any $w \in S_n(a)$, the Schubert cell $X^{\circ}_w(E\bull)$ is the $B$-orbit of $w\tilE$, and $\Fla$ is the disjoint union of these $B$-orbits as $w$ ranges over all of $S_n(a)$.  The proof of the following lemma can be found in ~\cite[Lemma 1]{BK06}.

\begin{lemma}The Schubert cell $X^{\circ}_w(F\bull)$ contains $\tilE$ if and only if $F\bull=pw^{-1}E\bull$ for some $p\in P$.\end{lemma}

Hence by Proposition \ref{Prop Ful_Mac}, $\prod_{k=1}^s [X_{w^k}]\neq0$ if and only if for generic $(p_1,p_2,\ldots,p_s)\in P^s$, the intersection $\bigcap_{k=1}^s X_{w^k}(p_k(w^k)^{-1}E\bull)$ is transverse at $\tilE$.  This fact gives motivation for the following definition.

\begin{definition}
The s-tuple $(w^1, w^2, \ldots, w^s)\in S_n(a)^s$ is \textbf{Levi-movable} or \textbf{L-movable} if
\begin{equation}\label{L_mov_dim}\sum_{k=1}^s\codim\, X_{w^k} =\dim\Fla\end{equation}

and for generic $(l_1, l_2, \ldots, l_s)\in L^s$, the intersection $\bigcap_{k=1}^s X_{w^k}(l_k(w^k)^{-1}E\bull)$ is transverse at the point $\tilE$.\end{definition}

The definition of $L$-movable can be defined in the same way for any flag variety $G/P$ and is used in this generality in Belkale-Kumar~\cite{BK06}.  One important property we will need is the following proposition:

\begin{prop}\label{Prop_MaxSame}
Consider the Grassmannian $\Gr(r,n)$.  If the s-tuple $(w^1,w^2,\ldots,w^s)\in S_n(r)^s$ satisfies equation \eqref{L_mov_dim} and $\prod_{k=1}^s [X_{w^k}]\neq0$ in $H^*(\Gr(r,n))$, then $(w^1,w^2,\ldots,w^s)$ is $L$-movable.\end{prop}

For proof see ~\cite[Lemma 19]{BK06}.


\section{Tangent spaces of flag varieties}\label{SECTION_TAN}

\quad Consider the injective map from $\Fla$ to the product $\prod_{i=1}^r \Gr(a_i,n)$ given by $$ V\bull \mapsto V_1 \times V_2 \times \cdots \times V_r.$$  It is well known the the tangent space of $\Gr(a_i,n)$ at the point $V_i$ is canonically isomorphic with $\Hom(V_i, \mathds{C}^n/V_i)$.  This induces an injective map on tangent spaces $$\Tv \Fla \hookrightarrow \bigoplus_{i=1}^r \Hom(V_i, \mathds{C}^n/V_i).$$ Let $i_j: V_j \hookrightarrow V_{j+1}$ and $\rho_j: \mathds{C}^n/V_j\twoheadrightarrow \mathds{C}^n/V_{j+1}$ denote the naturally induced maps from the flag structure of $V\bull$ and let $\phi := (\phi_1,\phi_2,\cdots,\phi_r)$ denote an element of $\bigoplus_{i=1}^r \Hom(V_i, \mathds{C}^n/V_i)$.  The tangent space of $\Fla$ at the point $V\bull$ is given by

\begin{equation}\label{Tangent}\Tv\Fla\simeq \{\phi \in \bigoplus_{i=1}^r \Hom(V_i,\mathds{C}^n/V_i)\ |\ \rho_j \circ \phi_j = \phi_{j+1} \circ i_j \ \forall j \in [r-1]\}.\end{equation}  If we consider the splitting \eqref{Q_splitting} together with the commuting conditions in \eqref{Tangent}, we have the following simplification of the tangent space of $\Fla$ at the point $\tilE$:

\begin{equation}\label{levitangent}
T_{\tilE}\Fla\simeq\bigoplus_{i=1}^{r}\Hom(Q_i,Q_{i+1}\oplus Q_{i+2} \oplus \cdots \oplus Q_{r+1})\simeq\bigoplus_{i<j}^{r+1} \Hom(Q_i, Q_j).\end{equation}

Note that $\Hom(Q_i, Q_j)$ is canonically isomorphic to the tangent space $T_{Q_i}(\Gr(b_i,Q_i\oplus Q_j))$.  We now describe the tangent space at a generic point of a Schubert cell $X_{w}^{\circ}(F\bull) \subseteq \Fla$.  To do this we need the notion of induced flags.  For any complete flag $F\bull\in\Fln$ and any subspace $V\subseteq \mathds{C}^n$, we have the induced complete flags on $V$ and $\mathds{C}^n/V$ given by the intersection of $V$ with $F\bull$ and the projection $\mathds{C}^n \twoheadrightarrow\mathds{C}^n/V$ of $F\bull$. We denote these induced flags by $F\bull(V)$ and $F\bull(\mathds{C}^n/V)$ respectively. We first recall a description of the tangent spaces of Schubert cells in Grassmannians.  For the proof see~\cite[Section 2.7]{So97}.

\begin{lemma}\label{Lemma_schu_tan_grass}
Let $r<n$ and $w\in S_n(r)$.  Let $F\bull\in \Fln$.  For any $V\in X_w^{\circ} (F\bull)\subseteq\Gr(r,n)$, the tangent space of the Schubert cell at the point $V$ is given by $$T_V X_w^{\circ}(F\bull)=\{\phi \in \Hom(V,\mathds{C}^n/V)\ |\ \phi(F_j(V))\subseteq F_{w(j)-j}(\mathds{C}^n/V) \quad \forall\ j\in [r] \}.$$\end{lemma}

We generalize this description to Schubert cells on the partial flag variety $\Fla$.  Let $P_{i}$ be the maximal parabolic identified with the set $\{a_i\}\subseteq [n-1]$.  For any $i \in [r]$, we have $W_P \subseteq W_{P_i}$ and the induced surjection $S_n(a)\twoheadrightarrow S_n(a_i)$. For any $w \in S_n(a)$, define $w_i \in S_n(a_i)$ to be the image of $w$ under this map.  Let $F\bull\in\Fln$ and let $V\bull \in X_w^{\circ} (F\bull)\subseteq \Fla$. The tangent space of the Schubert cell $X_{w}^{\circ} (F\bull)$ at the point $V\bull$ is given by $$\Tv X_w^{\circ}(F\bull)= \{\phi\in\Tv\Fla\ |\ \phi_i(F_j(V_i))\subseteq F_{w_{i}(j)-j}(\mathds{C}^n/V_i) \quad \forall\ i, j \}.$$ Consider the point $\tilE\in\Fla$ and choose $F\bull$ such that $\tilE\in X_w^{\circ}(F\bull)$.  We find that if $F\bull=lw^{-1}E\bull$ for some $l\in L$, then the space $T_{\tilE}X_w^{\circ}(F\bull)$ decomposes into a direct sum with respect to the decomposition \eqref{levitangent}.  Lemmas \ref{Lemma_VS_decomp} and \ref{trans_intersection} are basic facts about complex vector spaces and we leave the proofs as exercises.

\begin{lemma}\label{Lemma_VS_decomp}
Let $X=\bigoplus_{i=1}^rX_i$ be a complex vector space. Let $t=\{t_1,t_2,\ldots,t_r\}$ be distinct positive integers. Define the action of $t$ on $X$ by $t(\sum_{i=1}^r x_i)= \sum_{i=1}^r t_ix_i$.  Let $V$ be a subspace of $X$.  Then $t(V)=V$ if and only if $V=\bigoplus_{i=1}^r (V\cap X_i)$.\end{lemma}

\begin{lemma}\label{trans_intersection}
Let $X=\bigoplus_{i=1}^rX_i$ be a complex vector space. Let $Y^1, Y^2, \ldots, Y^s$ be vector subspaces such that $Y^k=\bigoplus_{i=1}^rY^k_i$ with $Y^k_i\subseteq X_i$ for all $i\in[r]$ and $k\in[s]$. If $$\dim\Big(\bigcap_{k=1}^s Y^k\Big)=\dim(X)-\sum_{k=1}^s\codim(Y^k)=0,$$ then $$\dim\Big(\bigcap_{k=1}^s Y^k_i\Big)=\dim(X_i)-\sum_{k=1}^s\codim(Y^k_i)=0$$ for each $i\in[s]$.\end{lemma}

\begin{definition}
For any partial flag variety $\Fla$ and $w\in S_n(a)$, define $$\Fl_{L,w}(n):=\{F\bull\in \Fln\ |\ F\bull= lw^{-1}E\bull\ \mbox{for some}\ l\in L\}$$\end{definition}

\begin{prop}\label{Prop_Tan_decomp}
Identify $$T_{\tilE}\Fla\simeq\bigoplus_{i<j}^{r+1}\Hom(Q_i, Q_j)$$ with respect to the splitting \eqref{Levi_split} of $\mathds{C}^n$ and for any $\{i<j\}\subseteq [r+1]$, identify
$$T_{Q_i} \Gr(b_i,Q_i\oplus Q_j)\simeq \Hom(Q_i, Q_j).$$ Let $w\in S_n(a)$ and $F\bull \in F\ell_{L,w}(n)$.  Then the following are true:

\begin{enumerate}
\item[(i)] $\displaystyle T_{\tilE}X_w^{\circ}(F\bull)= \bigoplus_{i<j}^{r+1} \big(\Hom(Q_i,Q_j)\cap T_{\tilE}X_w^{\circ}(F\bull)\big).$
\item[(ii)]$\Hom(Q_i,Q_j)\cap T_{\tilE}X_w^{\circ}(F\bull) = T_{Q_i}X_{w_{i,j}}^{\circ}(F\bull(Q_i\oplus Q_j))$ where $F\bull(Q_i\oplus Q_j)$ is the complete flag on $Q_i\oplus Q_j$ induced from $F\bull$.
\end{enumerate}\end{prop}

\begin{proof}Let $t=\{t_1,t_2,\ldots,t_{r+1}\}$ be a set of distinct positive integers. Let $t$ act on $\mathds{C}^n$ by $$t\big(\sum_{i=1}^{r+1}q_i\big)=\sum_{i=1}^{r+1}t_iq_i$$ with respect to the splitting \eqref{Q_splitting}. Since $F\bull=lw^{-1}E\bull$ for some $l\in L$, each $F_j(E_{a_i})$ is fixed by the $t$ action.  By Lemma \ref{Lemma_VS_decomp}, we have

$$F_j(E_{a_i})= \bigoplus_{m=1}^i (F_j(E_{a_i})\cap Q_m)$$ and
$$F_{w_i(j)-j}(\mathds{C}^n/E_{a_i}) = \bigoplus_{m=i+1}^{r+1} (F_{w_i(j)-j}(\mathds{C}^n/E_{a_i})\cap Q_m)$$

Therefore, for any $\phi=(\phi_1,\phi_2,\ldots,\phi_r)\in T_{\tilE}X_w^{\circ}(F\bull)$, the map $\phi_i$ can be written as the sum

\begin{equation}\label{Phi_QiQj}
\phi_i=\sum_{\substack{0<m_1\leq i\\i<m_2\leq r+1}}\phi_{m_1,m_2}\quad \mbox{where}\ \phi_{m_1,m_2}(F_{j}\cap Q_{m_1}) \subseteq F_{w_i(j)-j}\cap Q_{m_2}.\end{equation}

Note that $\phi_{m_1,m_2}\in \Hom(Q_{m_1},Q_{m_2})$.  Let $t'=\{t_{i,j}\}_{i<j}$ be a set of distinct positive integers and let $t'$ act on $T_{\tilE}\Fla$ by
$$t'\big(\sum_{i<j}^{r+1}\phi_{i,j}\big)=\sum_{i<j}^{r+1} t_{i,j}\phi_{i,j}$$ under the direct sum given in \eqref{levitangent}.  By equation \eqref{Phi_QiQj}, we have that $t'\big(T_{\tilE}X_w^{\circ}(F\bull)\big)=T_{\tilE}X_w^{\circ}(F\bull)$. Thus by Lemma \ref{Lemma_VS_decomp}, part $(i)$ of the proposition is proved.

\bigskip

\quad To prove part $(ii)$ of the proposition, we note that both sides of the equation are subspaces of $\Hom(Q_i,Q_j)$ of the same dimension and hence it suffices to show the $LHS\subseteq RHS$. For any $(\phi_1,\phi_2,\cdots,\phi_r)\in T_{\tilE}\Fla$ write $$(\phi_1, \phi_2, \cdots, \phi_r)=\sum_{i<j}^{r+1}\phi_{i,j}$$ under the decomposition \eqref{levitangent}. Fix $\{i<j\}\subseteq [r+1]$ and choose $m$ such that $i\leq m<j$.  Then we have

\begin{equation}\label{phi_m}\phi_m=\sum_{\substack{0<m_1\leq m\\m<m_2\leq r+1}}\phi_{m_1,m_2}.\end{equation}

Observe that $\phi_{i,j}$ is included in the above summation \eqref{phi_m}.  If $(\phi_1,\phi_2,\cdots,\phi_r)\in T_{\tilE}X_w^{\circ}(F\bull)$, then $$\phi_m(F_l(Q_1\oplus Q_2 \oplus \cdots \oplus Q_m))\subseteq F_{w_m(l)-l}(Q_{m+1}\oplus Q_{m+2}\oplus Q_{r+1}) \quad \forall\ l \in [a_m].$$  Define the set $M:=\{w(1), w(2),\ldots,w(a_{m})\}$ and for any
$k\in[b_i]$, let $$p_k:=\#\{\alpha \in M\ |\ \alpha \leq w(a_{i-1}+k) \}.$$  Note that $p_k$ is the smallest number such that the flag $F_{p_k}(Q_1\oplus Q_2 \oplus \cdots \oplus Q_m)$ induces the flag $F_k(Q_i)$ on $Q_i$ and that $w_m(p_k)=w(a_i+k)$.  By a basic calculation, one can show the flag $F_{w(a_i+k)-p_k}(Q_{m+1}\oplus Q_{m+2}\oplus \cdots \oplus Q_{r+1})$ induces the flag $F_l(Q_j)$ on $Q_j$ if and only if $w_{i,j}(k)=l+k$ . Hence the map $\phi_{i,j}$ in the decomposition \eqref{phi_m} of $\phi_m$ satisfies

\begin{equation}\label{wij_inclusion}\phi_{i,j}(F_k(Q_i))\subseteq F_{w_{i,j}(k)-k}(Q_j)\end{equation}

Note that this result is independent of choice of $m\in \{i,i+1,\ldots, j-1\}$. Since $\phi_{i,j}\in\Hom(Q_i,Q_j)$ that satisfy $\eqref{wij_inclusion}$ are exactly $\phi_{i,j}\in T_{Q_i}X^{\circ}_{w_{i,j}}(F\bull(Q_i\oplus Q_j))$, part $(ii)$ of proposition is proved.\end{proof}

\subsection*{Proof of $(ii)\Rightarrow (iii)$ in Theorem \ref{Thm1}:}

If $(w^1, w^2, \ldots, w^s)$ is $L$-movable, then for a generic point $(l_1,l_2, \ldots l_s)\in L^s$ we have $$\dim\Big(\bigcap_{k=1}^s T_{\tilE}X_{w^k}^{\circ}(F^k\bull)\Big)=\dim\Fla-\sum_{k=1}^s \codim(X_{w^k}^{\circ}(F^k\bull))=0$$ where $F^k\bull=l_k(w^k)^{-1}E\bull$.  By Proposition \ref{Prop_Tan_decomp}, $$\bigcap_{k=1}^s\big(T_{\tilE}X_{w^k}^{\circ}(F^k\bull)\big) =\bigoplus_{i<j}^{r+1}\Big(\bigcap_{k=1}^sT_{Q_i}X_{w^k_{i,j}}^{\circ}(F^k\bull(Q_i\oplus Q_j))\Big).$$ Hence by Lemma \ref{trans_intersection}, we have

\bigskip

$\displaystyle\dim\Big(\bigcap_{k=1}^sT_{Q_i}X_{w^k_{i,j}}^{\circ}(F^k\bull(Q_i\oplus Q_j))\Big)=$
$$\dim(\Hom(Q_i,Q_j))-\sum_{k=1}^s\codim\big(T_{Q_i}X_{w^k_{i,j}}^{\circ}(F^k\bull(Q_i\oplus Q_j))\big)=0$$ for any $\{i<j\}\subseteq [r+1]$.  This implies that the intersection $\bigcap_{k=1}^s X_{w^k_{i,j}}^{\circ}(F^k\bull(Q_i\oplus Q_j))$ is transverse at the point $Q_i\in\Gr(b_i,Q_i\oplus Q_j)$.  Choose $V\bull = Q_i\in \Gr(b_i,Q_i\oplus Q_j)$ and $\mathcal{F}=(F^1\bull(Q_i\oplus Q_j), F^2\bull(Q_i\oplus Q_j), \ldots, F^s\bull(Q_i\oplus Q_j))\in \Fl(Q_i\oplus Q_j)^s$ in Proposition \ref{Prop Ful_Mac}$(ii)$. This proves $(ii)\Rightarrow(iii)$ in Theorem \ref{Thm1}. \hfill$\Box$

\subsection*{Proof of $(iii)\Rightarrow (ii)$ in Theorem \ref{Thm1}:}
Assume that $\prod_{k=1}^s [X_{w^k_{i,j}}]$ is a nonzero multiple of a class of a point in $H^*(\Gr(b_i, b_i+b_j))$ for all $\{i<j\}\subseteq [r+1]$.  Since $\dim(Q_i\oplus Q_j)=b_i+b_j$, we identify the homogeneous space $\Gr(b_i, b_i+b_j)$ with $SL(Q_i\oplus Q_j)/P_{i,j}$.  Since $P_{i,j}$ is maximal, by Proposition \ref{Prop_MaxSame} we have that $(w^1_{i,j}, w^2_{i,j}, \ldots, w^s_{i,j})$ is $L_{i,j}$-movable. Hence for generic $\mathcal{F}_{i,j}\in \prod_{k=1}^{s} \Fl_{L_{i,j},w^k_{i,j}}(Q_i\oplus Q_j)$, we have
\begin{equation}\label{HomQQFG} \dim\Big(\bigcap_{k=1}^s T_{Q_i}X_{w^k_{i,j}}^{\circ}({F_{i,j}^k}\bull)\Big)= \dim(\Hom(Q_i,Q_j))-\sum_{k=1}^s\codim(X_{w^k_{i,j}}({F_{i,j}^k}\bull))=0\end{equation} for each $i<j$.  Let $U_{i,j}\subseteq\prod_{k=1}^{s} \Fl_{L_{i,j},w^k_{i,j}}(Q_i\oplus Q_j)$ denote an open set of flags which satisfy equation \eqref{HomQQFG}.  By Lemma \ref{Lemma_surj_flags} below, the map $\psi:= (\psi_{w^1},\psi_{w^2},\ldots,\psi_{w^s})$ is surjective and hence $\psi\circ\psi^{-1}(U_{i,j})=U_{i,j}$.  Since $\psi^{-1}(U_{i,j})$ is an open set, we can choose flags $(H\bull^1,\ldots,H\bull^s) \in \bigcap_{i<j}^{r+1}\psi^{-1}(U_{i,j})\subseteq\prod_{k=1}^s \Fl_{L,w^k}(n)$.  By Proposition \ref{Prop_Tan_decomp}, we have $$\dim\Big(\bigcap_{k=1}^s T_{\tilde{E\bull}}X_{w^k}^{\circ}(H^k\bull)\Big) = \sum_{i<j}^{r+1} \dim\Big( \bigcap_{k=1}^s T_{Q_i}X_{w^k_{i,j}}^{\circ}(H^k\bull(Q_i\oplus Q_j))\Big)=0.$$ Hence the intersection of the Schubert cells $X_{w^k}^{\circ}(H^k\bull)$ is transverse at the point $\tilde{E\bull}$.  Thus Proposition \ref{Prop Ful_Mac} proves $(iii)\Rightarrow (ii)$ in Theorem \ref{Thm1}. \hfill $\Box$

\bigskip

\begin{lemma}\label{Lemma_surj_flags} Fix $\{i<j\}\subseteq [r+1]$.  Let $L_{i,j}$ be the Levi subgroup of the parabolic $P_{i,j}
\subseteq SL(Q_i\oplus Q_j)$ which stabilizes the space $Q_i$.  Then for any $w\in
S_n(a)$, the map
$$\psi_w:\Fl_{L,w}(n) \rightarrow \Fl_{L_{i,j},w_{i,j}}(Q_i\oplus Q_j)$$
given by $\psi_w(F\bull)=F\bull(Q_i\oplus Q_j)$ is well defined and surjective.\end{lemma}

\begin{proof}Let $l=(g_1,g_2,\ldots,g_{r+1})\in L$ with respect to equation \eqref{Levi_split}.
By definition of $w_{i,j}$, we have $\psi_w(lw^{-1}E\bull)=(g_i,g_j)w_{i,j}^{-1}E\bull(Q_i\oplus
Q_j)$. Thus $\psi_w$ is well defined.  Since $\psi_w$ is $L$-equivariant and $L$ acts transitively
on $\Fl_{L_{i,j},w_{i,j}}(Q_i\oplus Q_j)$, we also have that $\psi_w$ is surjective. \end{proof}



\section{Structure coefficients}\label{SECTION_COEFF}

\quad Kleiman's transversality in ~\cite{Kl74} says that for a generic $s$-tuple of flags $\mathcal{F}=(F\bull^1,\ldots,F\bull^s)\in \Fln^s$, the intersection $\bigcap_{k=1}^sX^{\circ}_{w^k}(F^k\bull)$ is transversal and that it is dense in the intersection of projective varieties $\bigcap_{k=1}^sX_{w^k}(F^k\bull)$. If $(w^1,w^2,\ldots,w^s)\in S_n(a)^s$ such that $\prod_{k=1}^s[X_{w^k}]=c[pt]$ for some positive integer $c$, then for generic $\mathcal{F}\in \Fln^s$, we have

$$\left|\bigcap_{k=1}^sX^{\circ}_{w^k}(F^k\bull)\right|=\left|\bigcap_{k=1}^sX_{w^k}(F^k\bull)\right|=c.$$

The goal of this section is to produce a formula to write the structure coefficients $c$ corresponding to $L$-movable $s$-tuples as a product structure coefficients coming from Grassmannians.  To do this we consider the projection $f:\Fla\twoheadrightarrow \Gr(a_1,n)$.  We show that the number of points in an $L$-movable intersection of Schubert cells is equal to the number of points in the projected intersection in $\Gr(a_1,n)$ times the number of points in the fiber over any point the projection.  The techniques used in this section are inspired by techniques used by Belkale in his proof of Horn's conjecture in ~\cite{Be06}.
\smallskip

\quad Recall Definition \ref{def_induced_Weyl} of induced permutations.  Define the subset $\gamma:=\{a_1+1<a_1+2<\cdots<n\}\subseteq [n]$ and for any $w\in S_n(a)$, consider the induced permutation $w_{\gamma}\in S_{n-a_1}$.  For any point $V\in\Gr(a_1,n)$, the fiber $f^{-1}(V)$ is isomorphic to $\Fl(a_\gamma,n-a_1)$ where $a_\gamma=\{a_2-a_1<a_3-a_1<\cdots<a_r-a_1\}$.  It is easy to see that for any $w\in S_n(a)$ and $F\bull\in\Fln$ and $V\in f(X^{\circ}_{w}(F\bull))=X^{\circ}_{w_1}(F\bull)$, we have $$X^{\circ}_{w}(F\bull)\cap f^{-1}(V)\simeq X^{\circ}_{w_\gamma}(F\bull(\mathds{C}^n/V)).$$ The following lemma is immediate.

\begin{lemma}\label{Fiber_intersection} Let $(w^1,w^2,\ldots,w^s)\in S_n(a)^s$ and let $\mathcal{F}\in \Fln^s$ be such that $\bigcap_{k=1}^sX^{\circ}_{w_{1}^k}(F^k\bull)$ is not empty. For any $V\in\bigcap_{k=1}^sX^{\circ}_{w_{1}^k}(F^k\bull)$, we have
\begin{equation}\label{fiber_lemma_eq}\bigcap_{k=1}^sX^{\circ}_{w^k}(F^k\bull)\cap f^{-1}(V)\simeq \bigcap_{k=1}^sX^{\circ}_{w_{\gamma}^k}(F^k\bull(\mathds{C}^n/V)).\end{equation}\end{lemma}

Note that equation \eqref{fiber_lemma_eq} could possibly be empty.  However, we will show later in Proposition \ref{Prop_good_open}, that if $(w^1,w^2,\ldots,w^s)$ is $L$-movable, then we can choose $\mathcal{F}\in \Fln^s$ ``generic" enough so that \eqref{fiber_lemma_eq} is nonempty for all $V\in\bigcap_{k=1}^sX^{\circ}_{w_{1}^k}(F^k\bull)$.  We first show that an $L$-movable $s$-tuple induces a Levi-movable $s$-tuple on the fiber of $f$.  This will allow us to inductively apply Theorem \ref{Thm_constants} to $\Fl(a_\gamma,n-a_1)$ and hence reduce computing $L$-movable structure constants to the case of the Grassmannian.  We have the following relationships between the lengths of $w, w_\gamma,$ and $w_1$.

\begin{rem}For any $w\in S_n(a)$, we have
\begin{equation}\label{relation_global}\ell(w_\gamma)=\ell(w)-\ell(w_1)\end{equation}and
\begin{equation}\label{relation_induced}\ell((w_\gamma)_i)=\ell(w_{i+1})-\ell(w_1)+\sum_{k=1}^i \ell(w_{1,k}).\end{equation} for any $i\in[r-1]$\end{rem}

For any $V\in \Gr(a_1,n)$ and $w\in S_n(a_1)$, define $$Y_V^w:=\{F\bull \in \Fln\ |\ V\in X^{\circ}_w(F\bull)\}.$$

\begin{lemma}\label{Lemma_surj_flags2} For any $V\in \Gr(a_1,n)$ and $w\in S_n(a_1)$, the map $Y_V^w\rightarrow {\rm{F}}\ell(\mathds{C}^n/V)$ given by $F\bull \mapsto F\bull(\mathds{C}^n/V)$ is surjective.\end{lemma}

\begin{proof}Let $J=\{g\in \SLn\ |\ gV=V\}$.  It is easy to see that the map above is $J$-equivariant. Since $J$ acts transitively on ${\rm{F}}\ell(\mathds{C}^n/V)$, the map is surjective. \end{proof}

Let $P_1$ denote the maximal parabolic associated to $\Gr(a_1,n)$ and let $L_1$ be its Levi-subgroup. Let $R=L_1\cap P$ denote the parabolic subgroup of $L_1$.  We can identify $\Fl(a_{\gamma},n-a_1)$ with the homogeneous space $L_1/R$.  Let $L_R$ denote the Levi-subgroup of $R$.

\begin{prop}\label{L_induced}
If $(w^1,w^2,\ldots,w^s)$ is $L$-movable, then the following are true:
\begin{enumerate}
\item[(i)] The $s$-tuple $(w_1^1,w_1^2,\ldots,w_1^s)$ is $L_1$-movable.
\item[(ii)] The $s$-tuple $(w_\gamma^1,w_\gamma^2,\ldots,w_\gamma^s)$ is $L_R$-movable.\end{enumerate}\end{prop}

\begin{proof}We first show that for a generic choice of $\mathcal{F}\in \Fln^s$, the intersection $\bigcap_{k=1}^sX^{\circ}_{w_{1}^k}(F^k\bull)$ nonempty, finite and transverse.  Since $(w^1,w^2,\ldots,w^s)$ is $L$-movable, we can choose $\mathcal{F}$ so that this intersection is nonempty.  By Kleiman's transversality, the flags $\mathcal{F}$ can also be chosen so that the intersection is transverse.  Finally, by the numerical conditions \eqref{numerical_cond_eq}, the expected dimension of $\bigcap_{k=1}^sX^{\circ}_{w_{1}^k}(F^k\bull)$ is 0.  Hence the intersection is finite. Since $\bigcap_{k=1}^sX^{\circ}_{w_{1}^k}(F^k\bull)$ is an intersection of Schubert cells in a Grassmannian, by Proposition \ref{Prop_MaxSame}, the $s$-tuple $(w_1^1,w_1^2,\ldots,w_1^s)$ is $L_1$-movable.  This proves part$(i)$.


\quad For part$(ii)$, fix $V\in \Gr(a_1,n)$ and consider $\prod_{k=1}^sY_V^{w^k_1}\subseteq \Fln^s$.  By part$(i)$ and the assumption, for generic $\mathcal{F}\in \prod_{k=1}^sY_V^{w^k_1}$, we have that the intersections $\bigcap_{k=1}^sX^{\circ}_{w_{1}^k}(F^k\bull)$ and $\bigcap_{k=1}^sX^{\circ}_{w^k}(F^k\bull)$ are nonempty and transversal.  Since $f(\bigcap_{k=1}^sX^{\circ}_{w^k}(F^k\bull))\subseteq\bigcap_{k=1}^sX^{\circ}_{w_{1}^k}(F^k\bull)$, we can further assume that there exists a $V\bull\in\bigcap_{k=1}^sX^{\circ}_{w^k}(F^k\bull)$ such that $f(V\bull)=V$. By Lemma \ref{Fiber_intersection}, $\bigcap_{k=1}^sX^{\circ}_{w_{\gamma}^k}(F^k\bull(\mathds{C}^n/V))$ is nonempty and finite and by Lemma \ref{Lemma_surj_flags2}, the induced flags $\mathcal{F}(\mathds{C}^n/V)$ are generic in ${\rm{F}}\ell(\mathds{C}^n/V)^s$. Hence the intersection $\bigcap_{k=1}^sX^{\circ}_{w_{\gamma}^k}(F^k\bull(\mathds{C}^n/V))$ is transverse.  By Proposition \ref{Prop Ful_Mac}, $\prod_{k=1}^s[X_{w_{\gamma}^k}]$ is a nonzero multiple of a class of a point in $H^*({\rm F\ell}(a_\gamma,n-a_1))$.  By Theorem \ref{Thm1}, it suffices to check that $(w_\gamma^1,w_\gamma^2,\ldots,w_\gamma^s)$ satisfy the numerical conditions for $L_R$-movability. Since $(w^1,w^2,\ldots,w^s)$ is $L$-movable, we have the following numerical conditions:

\begin{eqnarray}
\sum_{k=1}^s \big(a_i(n-a_i)-\ell(w^k_i)\big)&=&a_i(n-a_i)\label{num1}\\
\sum_{k=1}^s \big(a_1(a_i-a_{i-1})-\ell(w^k_{1,i})\big)&=&a_1(a_i-a_{i-1})\label{num2}.
\end{eqnarray}

for any $i\in[r]$.  For any $i\in\{2,3,\ldots,r\}$ rewrite the dimension of $\Gr(\tilde{a}_{i-1},n-a_1)$ as

\begin{equation}\label{identity}\dim(\Gr(\tilde{a}_{i-1},n-a_1))=a_i(n-a_i)-a_1(n-a_1)+\sum_{k=1}^i a_1(a_k-a_{k-1}) \end{equation}

Combining \eqref{relation_induced},\eqref{num1},\eqref{num2} and \eqref{identity} shows that $(w_\gamma^1,w_\gamma^2,\ldots,w_\gamma^s)$ satisfies the numerical conditions for $L_R$-movability in Theorem \ref{Thm1}$(i)$. \end{proof}

Fix $(w^1,w^2,\ldots,w^s)$ to be $L$-movable.  We now show that for generic $\mathcal{F}=(F\bull^1,\ldots,F\bull^s)\in\Fln^s$, the intersection $\bigcap_{k=1}^sX^{\circ}_{w_{\gamma}^k}(F^k\bull(\mathds{C}^n/V))$ is nonempty for every $V\in\bigcap_{k=1}^s X^{\circ}_{w^k_1}(F^k\bull)$. Define the subvariety $Y\subseteq
\Gr(a_1,n)\times \Fln^s$ by the following: $$Y:=\{(V,\mathcal{F})\ |\ V\in \bigcap_{k=1}^s X^{\circ}_{w^k_1}(F^k\bull)\}.$$  By ~\cite[Proposition 8.1]{Be06} the variety $Y$ is irreducible and smooth.

\begin{definition}For any $(V,\mathcal{F})\in Y$, we say that $(V,\mathcal{F})$ has property $P1$ if the intersection $\bigcap_{k=1}^sX^{\circ}_{w_{\gamma}^k}(F^k\bull(\mathds{C}^n/V))$ is transversal and $\bigcap_{k=1}^sX^{\circ}_{w_{\gamma}^k}(F^k\bull(\mathds{C}^n/V))=\bigcap_{k=1}^sX_{w_{\gamma}^k}(F^k\bull(\mathds{C}^n/V)).$\end{definition}

Note that if $(V,\mathcal{F})$ has property $P1$, then $\bigcap_{k=1}^sX^{\circ}_{w^k}(F^k\bull)\cap f^{-1}(V)$ is not empty by Lemma \ref{Fiber_intersection}.

\begin{prop}Property $P1$ is an open condition on $Y$.\end{prop}

\begin{proof}Consider $Y$ as a fiber bundle on $\Gr(a_1,n)$ with fiber $\prod_{k=1}^sY_V^{w^k_1}$ over any point $V\in\Gr(a_1,n)$. Let $Z$ be the quotient flag bundle on $\Gr(a_1,n)$ with fiber ${\rm{F}}\ell(\mathds{C}^n/V)^s$ over any point $V\in\Gr(a_1,n)$.  By Lemma \ref{Lemma_surj_flags}, the fiber bundle map $\eta:Y\twoheadrightarrow Z$ given by $\mathcal{F}\mapsto \mathcal{F}(\mathds{C}^n/V)$ is surjective.  Choose an open set $U\subseteq\Gr(a_1,n)$ such that fiber bundle $Z$ is trivial. Over the set $U$, choose a local trivialization $$Z|_U \simeq U\times {\rm{F}}\ell(\mathds{C}^{n-a_1})^s.$$ Since $(w_\gamma^1,w_\gamma^2,\ldots,w_\gamma^s)$ is $L_R$-movable, there exists an open subset $O'\subset\Fl(\mathds{C}^{n-a_1})^s$ such that for every $\mathcal{H}\in O'$, the intersection $\bigcap_{k=1}^sX^{\circ}_{w_{\gamma}^k}(H^k\bull)$ is transversal and $\bigcap_{k=1}^sX^{\circ}_{w_{\gamma}^k}(H^k\bull) =\bigcap_{k=1}^sX_{w_{\gamma}^k}(H^k\bull)$ and define $O:=\bigcup_{g\in SL_{n-a_1}} gO'$. Clearly $O$ is $SL_{n-a_1}$-invariant under the diagonal action on ${\rm{F}}\ell(\mathds{C}^{n-a_1})^s$.  Consider the fiber bundle $\eta^{-1}(O)$ over $U$.  Since $O$ is $SL_{n-a_1}$-invariant, $\eta^{-1}(O)$ is independent of choice of local trivialization. It is easy to see that $\eta^{-1}(O)$ is an open set of $Y$ and every $(V,\mathcal{F})\in \eta^{-1}(O)$ satisfies property $P1$.\end{proof}

\begin{prop}\label{Prop_good_open}
Let $\tilde{O}\subseteq Y$ be an open subset of $Y$ such that every point in $\tilde{O}$ has property $P1$.  Let $g:Y\twoheadrightarrow \Fln^s$ be the projection of $Y$ onto its second factor. For generic $\mathcal{F}\in \Fln^s$, the set $g^{-1}(\mathcal{F})\subseteq \tilde{O}$.\end{prop}

\begin{proof}The fiber of $g$ over any point $\mathcal{F}$ is isomorphic to $\bigcap_{k=1}^s X^{\circ}_{w^k_1}(F^k\bull)$. Choose an open subset of $U_1\subseteq\Fln^s$ such that for every $\mathcal{F}\in U_1$, the set $g^{-1}(\mathcal{F})$ is finite.  Let $\tilde Y$ be the closure of $g(Y\backslash\tilde{O})$ in $\Fln^s$. Since $Y$ is irreducible, we have $\dim(Y\backslash\tilde{O})\geq \dim(\tilde{Y})$.  Since $g$ is generically finite to one, we have that $\dim(\Fln^s)>\dim(\tilde{Y})$.  Let $U_2$ be an open set in $\Fln^s\backslash\tilde Y$. For any $\mathcal{F} \in U_1\cap U_2$, we have $g^{-1}(\mathcal{F})\subseteq\tilde O$.\end{proof}

\subsection*{Proof of Theorem \ref{Thm_constants}:}

Chose $\mathcal{F}\in \Fln^s$ generically so that, $$\left|\bigcap_{k=1}^sX^{\circ}_{w_1^k}(F^k\bull)\right|=c_1\quad\mbox{and}\quad \left|\bigcap_{k=1}^sX^{\circ}_{w^k}(F^k\bull)\right|=c.$$ By Proposition \ref{Prop_good_open}, the flags $\mathcal{F}$ can also be generically chosen so that for any $V\in \bigcap_{k=1}^sX^{\circ}_{w_1^k}(F^k\bull)$ the point $(V,\mathcal{F})$ satisfies property $P1$. Therefore the map
$$f:\bigcap_{k=1}^sX^{\circ}_{w^k}(F^k\bull)\twoheadrightarrow \bigcap_{k=1}^sX^{\circ}_{w_1^k}(F^k\bull)$$
is surjective.  Since the number of points in each fiber $f^{-1}(V)$ is equal to $c_\gamma$, we have that $c=c_1\cdot c_\gamma.$\hfill $\Box$

\bigskip

Let $w_0$ be the longest element in $W=S_n$ and $w_a$ be the longest element in $W_P$. For any $w\in S_n(a)$, define $w^{\vee}:=w_0ww_a$.  Note that $w^{\vee}\in S_n(a)$.

\begin{cor}If $(w,u,v^{\vee})$ is $L$-movable, then $$c_{w,u}^v= c_{w_1,u_1}^{v_1}\cdot c_{w_{\gamma},u_{\gamma}}^{v_{\gamma}}.$$\end{cor}

\begin{proof}  By ~\cite[Lemma 16(d)]{BK06} the Poincar\'{e} pair $(w,w^{\vee})$ is $L$-movable.  By Proposition \ref{L_induced} we have that $(w_1,(w^{\vee})_1)$ and $(w_{\gamma},(w^{\vee})_{\gamma})$ are Levi-movable.  Hence $(w_1)^{\vee} = (w^{\vee})_1$ and $(w_{\gamma})^{\vee}=(w^{\vee})_\gamma$.  Apply Theorem \ref{Thm_constants} to the triple $(w,u,v^{\vee})$.\end{proof}

Recall that Proposition \ref{L_induced} says that if $(w,u,v^{\vee})$ is $L$-movable then $(w_\gamma,u_\gamma,v^{\vee}_\gamma)$ is $L_R$-movable.  Hence we can apply Theorem
\ref{Thm_constants} to $(w_\gamma,u_\gamma,v^{\vee}_\gamma)$.  This process gives an inductive way to write $c_{w,u}^v$ as a product of Littlewood-Richardson coefficients coming from the Grassmannians $\Gr(b_i,n-a_{i-1})$ where $i\in[r]$.

\begin{rem}Analogues of Theorem \ref{Thm_constants} exist for any projection $f_i:\Fla\twoheadrightarrow\Gr(a_i,n)$ and fiber $$f_i^{-1}(V)\simeq \Fl((a_1,\ldots,a_{i-1}),a_i)\times
\Fl((a_{i+1}-a_i,\ldots,a_r-a_i),n-a_i).$$ with corresponding induced coefficients.  The proof is similar to that of Theorem \ref{Thm_constants}.  Comparing these formulas gives many interesting relations between type A structure coefficients. \end{rem}

\subsection*{Acknowledgements}  The author is partially funded by NSF FRG grant DMS-05-54349.  I would like to thank my advisor Prakash Belkale for his patience and guidance in helping me write this paper.  I would also like to thank Richard Rimanyi for helpful suggestions in editing this paper.  I am grateful to the referees for their useful comments and remarks.  The results of this paper are part of on-going work for my Ph.D. thesis.

\end{document}